\documentclass[11pt]{article}
\usepackage{stmaryrd}
\usepackage{tipa}
\usepackage{amsfonts,amsmath,latexsym}
\usepackage{mathrsfs}   %不同于\mathcal or \mathfrak 之类的英文花体字体
\usepackage{amssymb}
\usepackage{CJK}
\usepackage{amsbsy}
\usepackage{indentfirst}

\def\a.s{\textsf{a. s}}

\def\and{\textsf{and}}

\newcommand{\beqnar}{\begin{eqnarray*}}
\newcommand{\eeqnar}{\end{eqnarray*}}
\newcommand{\ba}{\begin{array}}
\newcommand{\ea}{\end{array}}

\newcommand{\qed}{\nopagebreak\hspace*{\fill}
{\vrule width6pt height6ptdepth0pt}\par}
\allowdisplaybreaks[4]

\begin{document}
%首行缩进
%\setlength{\parindent}{2em}
%\begin{CJK*}{GBK}{song}
\title{\bf \Large On Convergence to Stochastic Integrals \footnote{Project supported by the National Natural Science
Foundation of China (No. 10871177), and the Specialized Research
Fund for the Doctor Program of Higher Education (No.
20090101110020).} }
 \date{\today}
\author{Zheng-Yan Lin, Han-Chao Wang\footnote{Corresponding author, hcwang06@gmail.com.}\\
Zhejiang University; Hangzhou, China, 310027}
\date{}
\maketitle $\mathbf{Abstract}$:  Weak convergence of various general
functionals of partial sums of dependent random variables to
stochastic integral now play a major role in the modern statistics
theory. In this paper, we obtain the weak convergence of various
general functionals of partial sums of causal process by means of
the method which was introduced in Jacod and Shiryaev (2003).

$\mathbf{Keyword}$: Weak convergence, causal process, stochastic
 integral.
%图书分类号
%{\bf2000  AMS-classification}: 60H30, 60H10, 60G46,
% 35C99, 58J65

 %Secondary : 60G15, 60G48.

%{\bf Actual version:} May 14th 2008

%\bigskip
%\bigskip

%\begin{itemize}
%作者介绍，脚标比较好？
%\end{itemize}
%这句貌似是换页？
%\vfill \eject

%花体章节分布
%\begin{center}
% {\large\bf\S\bf1\hspace*{3mm}  Introduction}
%\end{center}
\section{\bf{Introduction}}
\par $~~~$ Weak convergence of stochastic processes is a very
important and foundational theory in probability. In his classical
textbook, Billingsley (1968) gave a systematic theory of weak
convergence for stochastic processes.  In the theory, finite
dimensional distribution convergence and the tightness of stochastic
processes are crucial. Partial sum processes of random variables and
empirical processes are very important processes in the probability
and statistics. We can establish the weak convergence of  partial
sum processes to Brownian motion and empirical processes  to
Brownian bridge  by the classical method.

With the quick development of modern statistics and econometric
theory, the classical method has become difficult to deal with more
complex processes. For example,  sometimes  we need to prove a
convergence theorem about the stochastic integrals. However, we can
not get the weak convergence easily since it is difficult to compute
the finite dimensional distributions of the stochastic integrals.
The convergence theorem of stochastic integrals is  a core theory in
the unit root theory, which  is a hot topic in the econometric
theory. (c.f. Phillips (1987 a,b), (2007)). As mentioned in
Ibragimov and Phillips (2008), the results of this type can be used
in the study of transition behavior between regimes and marked
intervention policy. However, earlier authors only obtained some
results which describe the convergence to simple stochastic power
integrals by the classical weak convergence theory, since it is
complex to compute the finite dimensional distribution of more
general stochastic integral. (c.f. De Jong and Davidson (2000 a,
b)).

  In fact, when we intend to get a weak
convergence result for a  stochastic process sequence, we need to
complete two tasks : one is to prove the relative compactness of the
stochastic process sequence, and the other one  is to identify the
limiting process. In the Billingsley's classical method, tightness
is used to provide the relative compactness, while the convergence
of finite dimension distributions identifies the limiting process.
Jacod and Shiryaev (2003) developed a new approach to  weak
convergence of semimartingale sequences. Firstly, they introduced
three characteristics to replace the three terms: the drift,
variance of the Guassian part and L\'{e}vy measure, which
characterize the distribution of L\'{e}vy process. By means of these
three characteristics, one can show tightness of semimartingale
sequence. Secondly, they characterized the law of limiting process
as the unique solution of some martingale problem. In some special
cases, the unique solution of a martingale problem can be seen as a
unique solution of stochastic differential equations (for example,
the limiting process is a stochastic integral). Hence, Jacod and
Shiryaev's method is a powerful tool to prove the limit theorem
about the semimartingale. Because of some technical difficulty, this
method is rarely used in the statistics and econometric theory.
Ibragimov and Phillips (2008) used this method to obtain the weak
convergence of various general functionals of partial sums of linear
processes. This type of results can be used in the unit root theory,
where they deal with the functional of partial sum of linear
processes as a semimartingal and employ the Beveridge-Nelson
decomposition to compute  three predictable characteristics of the
underlying semimartingale.

In this paper, we extend these results to a causal process, which is
an important class of stationary processes. Wu (2005, 2007)
developed a complete method and theory about the causal process.  By
means of the martingale approximation developed by Wu (2005, 2007),
we obtain
 weak convergence of various general functionals of partial sums of
causal processes. In fact, the martingale approximation is an
extension of Beveridge-Nelson decomposition for a linear process. We
will use
 martingale approximation twice to obtain our results.

 The remainder of this paper is organized as follows.
        Section 2 gives a short introduction of the martingale convergence method developed by Jacod and Shiryaev (2003).
        Section 3 gives some definitions and notations about a
        causal process. Section 4 presents our main result.  The proof of the theorem will be
        collected in  Section 5. Section 6 gives a simple application
        of our result to unit root autoregression theory. Some
        discussion about the further research is given in Section 7.

\section{\bf{Martingale Convergence Method}}
   \subsection{\bf{ Definitions}}
     $~~~~~$In this subsection, we present
     some notations and preliminary results.

      Let $\mathbb{R}_{+}=[0,+\infty)$ and $\mathbb{Z}=\{\cdots, -2, -1, 0, 1, 2,
      \cdots\}$. $(\Omega, \mathscr{F}, \mathbb{F}=(\mathscr{F}_{t})_{t\ge 1},
      P)$ is a filtered probability space. $X$ is a semimartingale
      defined on $(\Omega, \mathscr{F}, \mathbb{F}=(\mathscr{F}_{t})_{t\ge 1},
      P)$.  Set $h(x)=x1_{|x|\le 1}$, and
           $$
\left\{
\begin{array}{ll}
\check{X}(h)_{t}=\sum_{s\leq t}[\Delta X_{s}-h(\Delta X_{s})],\\
X(h)=X-\check{X}(h),
\end{array}
\right. $$ $X(h)$ is a special semimartingale and we consider its
canonical decomposition:
  $$X(h)=X_{0}+M(h)+B(h),\eqno(2.1)$$
where $M(h)$ is its local martingale part, $B(h)$ is its finite
variation part.

$\bf{Definition~1} $  (Jacod and
     Shiryaev (2003) ) We call predictable characteristics of $X$
the triplet $(B,C,\nu)$ as follows:

  (1) $B$ is a predictable finite
variation process, namely the process $B=B(h)$ appearing in (2.1).

(2) $C=<X^{c},X^{c}>$ is a continuous process, where $X^{c}$ is the
continuous martingale part of $X$.

(3) $\nu$ is a predictable random measure on $\mathbb{R}_{+}\times
\mathbb{R}$, namely the compensator of the random measure $\mu^{X}$
associated to the jumps of $X$, $\mu^{X}$ is defined by
  $$\mu^{X}(\omega;dt,dx)=\sum_{s}1_{\{\triangle X_{s}(\omega)\neq 0\}}\varepsilon_{(s,\triangle X_{s})}(dt,dx),$$
where $\varepsilon_{a}$ denotes the Dirac measure at the point $a$.

 $\bf{Remark~1} $ By means of the truncation function $h(x)$,
the semimartingale $X$ can be divided into two parts: the jumps of
one part are greater than 1, and the jumps of the other's are not.
When a semimartingale's jumps are bounded, this semimartingale is a
special semimartingale, in the other words, it has unique canonical
decomposition, and hence  we can get an unique $B$ in Definition 1.
If the semimartingale is a special semimartingale, it is not
necessary to introduce the truncation function. In this paper, we
discuss such semimartingale, and so  do not introduce the truncation
function.

$\bf{Remark~2} $  The  predictable characteristics of semimartingale
$X$ are the counterpart of the drift, variance of Guassian part and
L\'{e}vy measure of independent increment process.  By means of
predictable characteristics, one can characterize the asymptotic
properties of the semimartingale.

If $\{Y_{k},~k\ge0\}$ is a discrete time semimartingale on
probability space $(\Omega, \mathscr{F}, P)$, we can write
    $$Y_{k}=\sum_{i=0}^{k}\eta_{i}=\eta_{0}+\sum_{i=1}^{k}m_{i}+\sum_{i=1}^{k}b_{i},$$
where $\eta_{0}=Y_{0}$, $\eta_{i}=Y_{i}-Y_{i-1}$,
$m_{i}=\eta_{i}-E(\eta_{i}|\mathscr{F}_{i-1})$ and
$b_{i}=E(\eta_{i}|\mathscr{F}_{i-1})$, $i\ge 1$.

Set $X_{s}=Y_{[s]}$. From  Definition 1, we can get the first and
second  predictable characteristic of $X_{s}$:
$$B_{s}=\sum_{i=0}^{[s]}b_{i},~C_{s}=\sum_{i=0}^{[s]}E(m_{i}^{2}|\mathscr{F}_{i-1}).\eqno(2.2)$$

The third predictable characteristic of $X_{s}$ is a compensated
random measure $\nu$. For a continuous function $g$ in $\mathbb{R}$,
we have
 $$\int_{0}^{s}\int_{\mathbb{R}}g(x)\nu(dx,dt)=\sum_{i=0}^{[s]}E(g(\eta_{i})|\mathscr{F}_{i-1}). \eqno(2.3)$$
 \subsection{\bf{Convergence of Semimartingales Using Predictable Characteristics}}

$\bf{Definition~2} $  (Jacod and
     Shiryaev (2003)) Let $X$ be a c\`{a}dl\`{a}g
     process and let $\mathcal {H}$ be the $\sigma-$field generated
     by $X(0)$ and $\mathcal {L}_{0}$ be the distribution of $X(0)$.
     A solution to the martingale problem associated with $(\mathcal {H},
     X)$ and $(\mathcal {L}_{0},B,C,\nu)$ (denoted by $\varsigma (\sigma(X_{0}),X| \mathcal {L}_{0}, B, C,
      \nu)$) is a probability measure
     $P$ on $(\Omega, \mathscr{F})$ such that $X$ is a
     semimartingale on $(\Omega, \mathscr{F}, P)$ with predictable
     characteristics $(B,C,\nu)$.

The limit process $X=(X(s))_{s\ge 0}$ appearing in this paper is the
canonical process $X(s,\alpha)=\alpha(s)$ for the element
$\alpha=(\alpha(s))_{s\ge 0}$ of $D(\mathbb{R}_{+})$. In  other
words, our limit process is defined on the canonical space
$(\mathbb{D}(\mathbb{R}_{+}), \mathscr{D}(\mathbb{R}_{+}),
\mathbf{D})$. For $a\ge 0$ and an element $(\alpha(s), s\ge 0)$ of
the Skorokhod space $\mathbb{D}(\mathbb{R}_{+})$, define
            $$S^{a}(\alpha)=\inf(s:|\alpha(s)|\ge a~\text{or}~|\alpha(s-)|\ge a),$$
              $$S_{n}^{a}=\inf(s:|X_{n}(s)|\ge a).$$
   In the paper, $\Rightarrow$ denotes convergence in distribution
   in an appropriate metric space, and $\xrightarrow[]{P}$ denotes
   convergence in probability. The following theorem gives
   sufficient conditions for the weak convergence of a sequence of
   square-integrable semimartingales. This theorem provides the
   basis for the study of asymptotic properties of functionals of
   partial sums.

   {\bf Theorem A} (Jacod and
     Shiryaev (2003)) Suppose that the following
     conditions hold:

     (i). The local strong majoration hypothesis: for all $a\ge 0$,
     there is an increasing continuous and deterministic function
     $F(a)$ such that the stopped processes $Var(B)^{S^{a}}$,
     $C^{S^{a}}$ and $(|x^{2}*\nu|)^{S^{a}}$ are strongly majorized
     by $F(a)$.

     (ii). The local condition on big jumps: for all $a\ge 0$, $t\ge
     0$,
           $$\lim_{b\uparrow \infty}\sup_{\alpha\in \Omega}|x^{2}|1_{\{|x|>b\}}*\nu_{t\wedge S^{a}}(\alpha)=0.$$

      (iii). Local uniqueness for the martingale problem $\varsigma (\sigma(X_{0}),X| \mathcal {L}_{0}, B, C,
      \nu)$; We denote by $Q$ the unique solution to this problem.

      (iv). Continuity condition: for all $t\in D$, $g\in
      \mathbb{C}(R)$, the function $\alpha\rightsquigarrow B_{t}(\alpha), C_{t}(\alpha),
      g*\nu_{t}(\alpha)$ are Skorokhod-continuous on
      $\mathbb{D}(R)$, where $D$ is a dense subset of $\mathbb{R}_{
      +}$.

      (v). $\mathcal {L}_{0}^{n}\rightarrow\mathcal {L}_{0}$ weakly
      as $n\rightarrow\infty$.

       (vi). $g*\nu_{t\wedge S_{n}^{a}}^{n}-(g*\nu_{t\wedge S^{a}})\circ(X^{n})\xrightarrow[]{P}
       0$ for all $t\in D$, $a>0$, $g\in \mathbb{C}_{+}(R)$;\\
         $~~~~~~~~~~~$ $\sup_{s\le t}|B_{t\wedge S_{n}^{a}}^{n}-(B_{t\wedge S^{a}})\circ(X^{n})|\xrightarrow[]{P}
         0$ for all $t>0$, $a>0$;\\
         $~~~~~~~~~~~$ $C_{t\wedge S_{n}^{a}}^{n}-(C_{t\wedge S^{a}})\circ(X^{n})\xrightarrow[]{P}
         0$ for all $t\in D$, $a>0$;\\
         $~~~~~~~~~~~$ $\lim _{b\uparrow\infty}\lim\sup_{n}P(|x^{2}|1_{\{|x|>b\}}*\nu_{t\wedge
         S_{n}^{a}}^{n}>\varepsilon)=
         0$ for all $t>0$, $a>0$, $\varepsilon>0$.\\Then $\mathscr{L}(X^{n})\Rightarrow Q$.

 \subsection{\bf{Uniqueness Conditions for Homogenous Diffusion Processes }}

$~~~~~$The limiting process in Theorem A  usually can be seen as a
homogenous diffusion process.  In this paper, the  stochastic
differential equation which is discussed has a homogenous diffusion
process  solution, we need a theorem to assure that this stochastic
differential equation has a unique and measurable solution.

Consider the stochastic differential equation:
$$
\begin{cases}
dX_{1}(t)=\lambda g(X_{2}(t))dt+\sigma f(X_{2}(t))dB(t),\\
dX_{2}(t)=dB(t).
\end{cases}
\eqno(2.4)$$  A solution to (2.4) is a two-dimensional
semimartingale $X:=(X_{1},X_{2})$ with the predictable
characteristics $B(X)$ and $C(X)$, where, for an element
$\alpha(s)=(\alpha_{1}(s),\alpha_{2}(s)) $ in
$\mathbb{D}(\mathbb{R}^{2})$,
\begin{eqnarray*}& &B(s,\alpha)=(\int_{0}^{s}g(\alpha_{2}(v))dv, 0),\\
                 & & C(s,\alpha)=\begin{bmatrix}
\int_{0}^{s}f^{2}(\alpha_{2}(v))dv&\int_{0}^{s}f(\alpha_{2}(v))dv\\\int_{0}^{s}f(\alpha_{2}(v))dv&s
\end{bmatrix}.\end{eqnarray*}

 {\bf Theorem B}
(Ibragimov and Phillips (2008) ) Suppose that

(i) The functions $f(x)$ and $g(x)$ are locally Lipschitz
continuous, i.e. for every $N\in \mathbf{N}$, there exists a
constant $K_{N}$ such that
$$ |f(x)-f(y)|\le K_{N}|x-y|,~~|g(x)-g(y)|\le K_{N}|x-y|$$
for all $|x|\le N$,  $|y|\le N$.

(ii) $f$ and $g$ satisfy the growth condition: there exists $K>0$
such that
    $$|f(x)|\le e^{K|x|},~~|g(x)|\le e^{K|x|}.$$
Then the stochastic differential equation (2.4) has a unique
solution. In  other words, the martingale problem $\varsigma
(\sigma(X_{0}),X| \mathcal {L}_{0}, B, C,
      0)$ has an unique solution.

\section{\bf{Causal Process and Martingale
Approximation}}
  We call $\{X_{n},n\ge 1\}$,  a causal process if $X_{n}$ has the form
  $$X_{n}=g(\cdots,\varepsilon_{n-1},\varepsilon_{n}),$$
where  $\{\varepsilon_{n}; n\in Z\}$ is mean zero, independent and
identically distributed  random variables and $g$ is a measurable
function. Causal process is a very
  important example of stationary process. It has been widely used
  in practice,  and contains many important statistical models, such as
ARCH models, threshold AR (TAR) and so on.  Asymptotic behavior of
the sums of causal processes, $S_{n}=\sum_{i=1}^{n}X_{i}$, are
important subjects in both practice and theory.

Recall that  $Z\in L^{p}$ $(p>0)$ if
       $||Z||_{p}=[E(|Z|^{p})]^{1/p}<\infty$ and write
       $||Z||=||Z||_{2}$.

To study the asymptotic property of the sums of causal processes,
martingale approximation is an effective method. Roughly speaking,
martingale approximation is to find a martingale $M_{n}$, such that
the error $\parallel S_{n}-M_{n}\parallel_{p}$ is small in some
sense. We list the notations  used in the following part:

$\bullet$
$\mathscr{F}_{k}=(\cdots,\varepsilon_{k-1},\varepsilon_{k})$.

$\bullet$ Projections $\mathcal
{P}_{k}Z=E(Z|\mathscr{F}_{k})-E(Z|\mathscr{F}_{k-1})$, $Z\in L^{1}$.

$\bullet$  $D_{k}=\sum_{i=k}^{\infty}{P}_{k}X_{i}$,
$M_{k}=\sum_{i=1}^{k}D_{i}$, $R_{k}=S_{k}-M_{k}$.

$\bullet$ $H_{k}=\sum_{i=1}^{\infty}E(X_{k+1}|\mathscr{F}_{k})$.

$\bullet$ $\theta_{n,p}=||\mathcal{P}_{0}X_{n}||_{p},$,
$\Lambda_{n,q}=\sum_{i=0}^{n}\theta_{i,q}$, $n> 0$. Let
$\theta_{n,p}=0=\Lambda_{n,p}$ if $n<0$.

$\bullet$ $\Theta_{m,p}=\sum_{i=m}^{\infty}\theta_{i,p}$ .

$\bullet$ $B$: standard Brownian motion.

 $M_{k}$ is a martingale, we will use $M_{k}$ to
approximate  sum $S_{k}$. Throughout the paper, we assume that
$D_{k}$ converges almost surely.

Linear process is a very important example of causal processes, and
many methods are developed to discuss it. Phillips and Solo (1992)
studied the Beveridge-Nelson decomposition of linear processes, and
then obtained some  asymptotic results. This method is used to
obtain the asymptotic results of short memory linear processes.

Suppose that $U_{n}$ is the linear process
$U_{n}=\sum_{i=0}^{\infty}a_{i}\varepsilon_{n-i}$.  Applying the
Beveridge-Nelson polynomial decomposition, one can get
      $$U_{n}=(\sum_{i=0}^{\infty}a_{i})\varepsilon_{n}+\widetilde{\varepsilon}_{n-1}-\widetilde{\varepsilon}_{n},\eqno(3.1)$$
where
$\widetilde{\varepsilon}_{n}=\sum_{i=0}^{\infty}\widetilde{a}_{i}\varepsilon_{n-i}$,
$\widetilde{a}_{i}=\sum_{k=i+1}^{\infty}a_{k}$.

In fact, the martingale approximation introduced above is an
extension of the Beveridge-Nelson decomposition to causal processes.
Applying the martingale approximate to $U_{n}$, we get
      $$\sum_{i=k}^{\infty}{P}_{k}U_{i}=(\sum_{i=0}^{\infty}a_{i})\varepsilon_{k},~~U_{k}-\sum_{i=k}^{\infty}{P}_{k}U_{i}=\widetilde{\varepsilon}_{k-1}-\widetilde{\varepsilon}_{k}.$$
From Wu (2005), we have
                     $$||\mathcal {P}_{0}U_{n}||_{p}=c_{0}|a_{n}|,~~~~~c_{0}=||\varepsilon_{0}-\varepsilon_{0}'||_{p},\eqno(3.2)$$
where $\varepsilon_{0}'$ is the independent copy of
$\varepsilon_{0}$.

\section{\bf{Main Result}}
 {\bf Assumption 1.} $X_{0}\in \mathcal {L}^{q}$, $q>2$, and
$\Theta_{n,q^{*}}=O(n^{1/q^{*}-1/2}(\log n)^{-1})$, where
$q^{*}=\min(q,4)$.

{\bf Assumption 2.}
          $$\sum_{k=1}^{\infty}||E(D_{k}^{2}|\mathscr{F}_{0})-\sigma^{2}||_{q^{*}/2}<\infty,$$
where $||D_{k}||=\sigma$.

 {\bf Assumption 3.}
$$\sum_{k=0}^{\infty}\sum_{i=1}^{\infty}||E(X_{k}X_{k+i}|\mathscr{F}_{0})-E(X_{k}X_{k+i}|\mathscr{F}_{-1})||_{4}<\infty,$$
and
$$||\sum_{k=0}^{\infty}\sum_{i=1}^{\infty}E(X_{k}X_{k+i}|\mathscr{F}_{0})||_{3}<\infty.$$
 {\bf Theorem 1} Let $f:\mathbb{R}\rightarrow \mathbb{R}$ be a twice
 continuously differentiable function such that $f'$ satisfies $|f'(x)|\le
 K(1+|x|^{\alpha})$ for some positive constants $K$
and $\alpha$ and all $x\in \mathbb{R}$. Suppose that $X_{t}$ is a
causal process  satisfies Assumptions 1$\thicksim$3. Then
       $$\frac{1}{\sqrt{n}}\sum_{t=2}^{[nr]}f(\frac{1}{\sqrt{n}}\sum_{i=1}^{t-1}X_{i})X_{t}\Rightarrow \lambda\int_{0}^{r}f'(B(v))dv+\sigma\int_{0}^{r}f(B(v))dB(v),\eqno(4.1)$$
where $\lambda=\sum_{j=1}^{\infty}EX_{0}X_{j}$.

{\bf Remark 1} The assumptions on the function $f$ in this paper is
the same as that of Ibragimov and Phillips (2008). Assumption 1 on
the causal process is more wild than Ibragimov and Phillips
(2008)'s. In their paper, they assume
$\sum_{i=1}^{\infty}i|a_{i}|<\infty$, while our condition is weaker
from (3.2).

{\bf Remark 2}  Assumption 3 is not strong as well. We have
              $$||\sum_{k=0}^{\infty}\sum_{i=1}^{\infty}E(X_{k}X_{k+i}|\mathscr{F}_{0})||_{3}\le O(\sum_{r=0}^{\infty}\sum_{k=0}^{\infty}\sum_{j=k+1}^{\infty}|a_{j}||\widetilde{a}_{j+r}|).$$
From Lemma E5 of Ibragimov and Phillips (2008),
$\sum_{i=1}^{\infty}i|a_{i}|<\infty$ implies
$\sum_{r=0}^{\infty}\sum_{k=0}^{\infty}\sum_{j=k+1}^{\infty}|a_{j}||\widetilde{a}_{j+r}|<\infty$.

\section{\bf{The Proof of Main Result}}

$~~~~~$Firstly, we introduce two lemmas.

{\bf Lemma 1} (Wu (2007) ) Assume that $E[X_{0}]=0$, $X_{0}\in
\mathcal {L}^{q}$, $q>1$, let $q'=\min(2,q)$, and
$\Theta_{0,q}=\sum_{i=0}^{\infty}\theta_{i,q}<\infty$, then
$$\parallel\max_{k\le n}|S_{k}|\parallel_{q}\le \frac{qB_{q}}{q-1}n^{1/q'}\Theta_{0,q},$$
where $B_{q}=18q^{3/2}(q-1)^{-1/2}$ if $q\in (1,2)\cup (2,\infty)$
and $B_{q}=1$ if $q=2$.

{\bf Lemma 2} (Wu (2007) )
    Under Assumptions 1 and 2, there exists a standard Brownian
    motion $B$ on a richer probability space such that
     $$|S_{n}-B(\sigma^{2}n)|=O_{a.s.}(n^{1/4}(\log n)^{1/2}(\log\log n)^{1/4}).$$

 Set
      $$X_{n}(s)=(\frac{1}{\sqrt{n}}\sum_{t=2}^{[ns]}f(\frac{1}{\sqrt{n}}\sum_{i=1}^{t-1}X_{i})X_{t}, \frac{1}{\sqrt{n}}\sum_{t=1}^{[ns]}X_{t})=:(X_{n}^{1}(s), X_{n}^{2}(s))$$
and
      $$X(s)=( \lambda\int_{0}^{s}f'(B(v))dv+\sigma\int_{0}^{s}f(B(v))dB(v), B(s)
      )=:(X^{1}(s),X^{2}(s)).$$  By (2.2) and (2.3), we can get the first two predictable characteristics of
      $X_{n}$ as follows:
      \begin{eqnarray*}& &B_{n}(s)=(\frac{1}{\sqrt{n}}\sum_{t=2}^{[ns]}f(\frac{1}{\sqrt{n}}\sum_{i=1}^{t-1}X_{i})(X_{t}-D_{t}),\frac{1}{\sqrt{n}}\sum_{t=1}^{[ns]}(X_{t}-D_{t})
),\\
    & &C_{n}(s)=\begin{bmatrix}
  C_{n}^{11}(s)&C_{n}^{12}(s)\\ C_{n}^{21}(s)&C_{n}^{22}(s)
\end{bmatrix},\\
& &
C_{n}^{11}(s)=\frac{1}{n}\sum_{t=2}^{[ns]}f^{2}(\frac{1}{\sqrt{n}}\sum_{i=1}^{t-1}X_{i})E(D_{t}^{2}|\mathscr{F}_{t-1}),\\
& &
C_{n}^{22}(s)=\frac{1}{n}\sum_{t=1}^{[ns]}E(D_{t}^{2}|\mathscr{F}_{t-1}),\\
& &
C_{n}^{12}(s)=C_{n}^{21}(s)=\frac{1}{n}\sum_{t=2}^{[ns]}f(\frac{1}{\sqrt{n}}\sum_{i=1}^{t-1}X_{i})E(D_{t}^{2}|\mathscr{F}_{t-1}).\end{eqnarray*}
 The process
$X(s)=(X^{1}(s),X^{2}(s))$,  an element of the Skorokhod space
$D(\mathbb{R}_{+})$,  is a solution to the stochastic differential
process
$$
\begin{cases}
dX^{1}(t)=\lambda f'(X^{2}(t))dt+\sigma f(X^{2}(t))dB(t),\\
dX^{2}(t)=dB(t).
\end{cases}
\eqno(5.1)$$

The predictable characteristics of $X$ are $(B(X),C(X),0)$:
\begin{eqnarray*}& &B(s,\alpha)=(\lambda\int_{0}^{s}f'(\sigma\alpha_{2}(v))dv,0),\\& &C(s,\alpha)=\begin{bmatrix}
\sigma^{2}\int_{0}^{s}f^{2}(\alpha_{2}(v))dv&\sigma\int_{0}^{s}f(\alpha_{2}(v))dv\\\sigma\int_{0}^{s}f(\alpha_{2}(v))dv&\sigma^{2}s
\end{bmatrix}.\end{eqnarray*}

To prove the main result, we need to verify the conditions in
Theorem A.

 By the similar argument in Ibragimov and Phillips (2008)
(3.15-3.17), we can get that condition (i) in Theorem A is
satisfied.

Since the limiting process is continuous, condition (ii) and
condition (iv) in Theorem A don't need to be verified.

 Under the assumptions of the
theorem, function $f(x)$ and $f'(x)$ are locally Lipschitz
continuous and satisfy growth condition. From Theorem B, the
stochastic differential equation has a unique solution. In other
words, the martingale problem $\varsigma (\sigma(X_{0}),X| \mathcal
{L}_{0}, B, C,
      \nu)$ have unique solution. We can get that condition (iii) in Theorem A is
satisfied.

Since $\mathcal {L}_{0}=0$, $\mathcal {L}^{n}_{0}=0$, it suffices to
check condition (vi) in Theorem A.

 From Jacod and Shiryaev (2003), if we can show
      $$\sup_{0<s\le N}|\Delta X_{n}(s)|\xrightarrow[]{P} 0~~\text{for all}~N\in \mathbf{N}, \eqno(5.2)$$
then $g*\nu_{t\wedge S_{a}^{n}}^{n}-(g*\nu_{t\wedge
S_{a}})\circ(X^{n})\xrightarrow[]{P}
       0$ for all $t\in D$, $a>0$, $g\in \mathbb{C}_{+}(R)$.

Secondly, we need to compute the terms of $B(s)\circ X_{n}$,
$C(s)\circ X_{n}$ :
  $$B(s)\circ X_{n}= \frac{\lambda}{n}\sum_{t=2}^{[ns]}f'(\frac{1}{\sqrt{n}}\sum_{i=1}^{t-1}X_{i})+\frac{\lambda}{n}f'(\frac{1}{\sqrt{n}}\sum_{i=1}^{[ns]}X_{i})(ns-[ns]),$$
  $$C^{11}(s)\circ X_{n}=\frac{\sigma^{2}}{n}\sum_{t=2}^{[ns]}f^{2}(\frac{1}{\sqrt{n}}\sum_{i=1}^{t-1}X_{i})+\frac{\sigma^{2}}{n}f^{2}(\frac{1}{\sqrt{n}}\sum_{i=1}^{[ns]}X_{i})(ns-[ns]),$$
    $$C^{12}(s)\circ X_{n}=C^{21}(s)\circ X_{n}=\frac{\sigma^{2}}{n}\sum_{t=2}^{[ns]}f(\frac{1}{\sqrt{n}}\sum_{i=1}^{t-1}X_{i})+\frac{\sigma^{2}}{n}f(\frac{1}{\sqrt{n}}\sum_{i=1}^{[ns]}X_{i})(ns-[ns]).$$
We divide the proof into three steps:

 {\bf Step 1} We  prove
 $$\sup_{0<s\le N}|C_{n}^{ij}(s)-C^{ij}(s)\circ X_{n}|\xrightarrow[]{P}0, \eqno(5.3)$$
where $N$ is a positive integer.  Since the proofs  for $i,j=1,2$
are similar, we only consider the case of $i=j=1$. In fact we need
to show that
      $$\sup_{0<s\le N}|\frac{1}{n}\sum_{t=2}^{[ns]}f^{2}(\frac{1}{\sqrt{n}}\sum_{i=1}^{t-1}X_{i})(E(D_{t}^{2}|\mathscr{F}_{t-1})-\sigma^{2})-\frac{\sigma^{2}}{n}f^{2}(\frac{1}{\sqrt{n}}\sum_{i=1}^{[ns]}X_{i})(ns-[ns])|\xrightarrow[]{P}0,\eqno(5.4)$$
Firstly, we prove
               $$\sup_{0<s\le N}|\frac{1}{n}\sum_{t=2}^{[ns]}f^{2}(\frac{1}{\sqrt{n}}\sum_{i=1}^{t-1}X_{i})\sigma^{2}-\frac{1}{n}\sum_{t=2}^{[ns]}f^{2}(\frac{1}{\sqrt{n}}\sum_{i=1}^{t-1}D_{i})\sigma^{2}|\xrightarrow[]{P}0.\eqno(5.5)$$
Since
$$\sup_{0<s\le N}|\frac{1}{n}\sum_{t=2}^{[ns]}f^{2}(\frac{1}{\sqrt{n}}\sum_{i=1}^{t-1}X_{i})-\frac{1}{n}\sum_{t=2}^{[ns]}f^{2}(\frac{1}{\sqrt{n}}\sum_{i=1}^{t-1}D_{i})|
                   \le \max_{1\le t\le
                   nN}|f^{2}(\frac{1}{\sqrt{n}}\sum_{i=1}^{t-1}X_{i})-f^{2}(\frac{1}{\sqrt{n}}\sum_{i=1}^{t-1}D_{i})|,$$
and $f$ is a uniform continuous function, we can get (5.5) by
                   $$\max_{1\le t\le
                   nN}|\frac{1}{\sqrt{n}}\sum_{i=1}^{t-1}X_{i}-\frac{1}{\sqrt{n}}\sum_{i=1}^{t-1}D_{i}|\xrightarrow[]{P}0.\eqno(5.6)$$
For any $\varepsilon>0$, by Lemma 1, for $2<q<4$, we have
            $$P(\frac{1}{\sqrt{n}}\max_{1\le t\le
                   nN}|\sum_{i=1}^{t-1}X_{i}-\sum_{i=1}^{t-1}D_{i}|>\varepsilon)\le \frac{E[\sum_{t=1}^{nN}(X_{t}-D_{t})]^{2}}{n\varepsilon^{2}}\le C\frac{n^{1/q'}(\log n)^{-1}}{n\varepsilon^{2}},$$
which implies (5.6), then we obtain (5.5).

By the martingale version of  the Skorokhod representation theorem,
on a richer probability space, there exist a standard Brownian
motion $B$ and nonnegative random variables
$\tau_{1},\tau_{2},\cdots$ with partial sums
$T_{k}=\sum_{i=1}^{k}\tau_{i}$ such that for $k\ge1$,
$T_{k}-k\sigma^{2}=o_{\text{a.s.}}(k^{2/q})$ and
      $M_{k}=B(T_{k})$,
      $E(\tau_{k}|\mathscr{F}_{k-1})=E(D_{k}^{2}|\mathscr{F}_{k-1})$.
      (cf. Wu (2007)).

For $\frac{T_{k-1}}{n}<s\le \frac{T_{k}}{n}$, we consider
$$\mathcal
{I}_{n}(s)=\frac{1}{n}\sum_{t=2}^{[ns]}f^{2}(\frac{1}{\sqrt{n}}\sum_{i=1}^{t-1}D_{i})\sigma^{2}+\frac{\sigma^{2}}{n}f^{2}(\frac{1}{\sqrt{n}}\sum_{i=1}^{[ns]}D_{i})(ns-[ns]).$$
By the  martingale version of  the Skorokhod representation theorem,
we have
   $$\mathcal
{I}_{n}(s)=\sum_{t=2}^{k-1}f^{2}(B(\frac{T_{t-1}}{n}))\frac{\sigma^{2}}{n}+\frac{\sigma^{2}}{n}f^{2}(B(\frac{T_{k-1}}{n}))(ns-[ns]).$$
Since $T_{k}-k\sigma^{2}=o_{\text{a.s.}}(k^{2/q})$,
    \begin{eqnarray*}\max_{t\le k}|B(\frac{T_{t}}{n})-B(\frac{\sigma^{2}t}{n})|&\le& \max_{t\le k}\sup_{|x-\sigma^{2}t|\le k^{2/q}}|B(\frac{x}{n})-B(\frac{\sigma^{2}t}{n})| \\
                                                                              &\le & o_{\text{a.s.}}(k^{1/q}\sqrt{\log k}).\end{eqnarray*}
By the uniform continuous property of $f$ and the similar argument
in (5.5), we have
$$\sup_{0<s\le N}|\mathcal
{I}_{n}(s)-\sum_{t=2}^{k-1}f^{2}(B(\frac{\sigma^{2}
(t-1)}{n}))\frac{\sigma^{2}}{n}+\frac{\sigma^{2}}{n}f^{2}(B(\frac{\sigma^{2}
(k-1)}{n}))(ns-[ns])|\xrightarrow[]{P}0.$$
 By the  Riemann approximation of stochastic integral, and the continuity of
Brownian motion's paths, we have
$$\sup_{0<s\le N}|\mathcal {I}_{n}(s)-\int_{0}^{s}f^{2}(B(v))dv|\xrightarrow[]{P}0.\eqno(5.7)$$
For
$\frac{1}{n}\sum_{t=2}^{[ns]}f^{2}(\frac{1}{\sqrt{n}}\sum_{i=1}^{t-1}X_{i})E(D_{t}^{2}|\mathscr{F}_{t-1})$,
and by noting $M_{k}=B(T_{k})$, we have
\small{ \begin{eqnarray*}& &\sup_{0<s\le N}| \frac{1}{n}\sum_{t=2}^{[ns]}f^{2}(\frac{1}{\sqrt{n}}\sum_{i=1}^{t-1}X_{i})E(D_{t}^{2}|\mathscr{F}_{t-1})- \frac{1}{n}\sum_{t=2}^{[ns]}f^{2}(\frac{1}{\sqrt{n}}B(T_{t-1}))E(D_{t}^{2}|\mathscr{F}_{t-1})|\\
                  &=& \sup_{0<s\le N}|\frac{1}{n}\sum_{t=2}^{[ns]}f^{2}(\frac{1}{\sqrt{n}}\sum_{i=1}^{t-1}X_{i})E(T_{t}-T_{t-1}|\mathscr{F}_{t-1})-
                  \frac{1}{n}\sum_{t=2}^{[ns]}f^{2}(\frac{1}{\sqrt{n}}B(T_{t-1}))E(T_{t}-T_{t-1}|\mathscr{F}_{t-1})|~~(5.8)\\
                  &\xrightarrow[]{P}&0. \end{eqnarray*}}
by Lemma 2.

By the  Riemann approximation of stochastic integral and the
Approximated Laplacians property (cf. Dellacherie and Meyer (1982)),
 we obtain
 $$\sup_{0<s\le N}|\frac{1}{n}\sum_{t=2}^{[ns]}f^{2}(\frac{1}{\sqrt{n}}B(T_{t-1}))E(T_{t}-T_{t-1}|\mathscr{F}_{t-1})-\int_{0}^{s}f^{2}(B(v))dv|\xrightarrow[]{P}0,$$
and combining with (5.8), we have
$$\sup_{0<s\le N}|\frac{1}{n}\sum_{t=2}^{[ns]}f^{2}(\frac{1}{\sqrt{n}}\sum_{i=1}^{t-1}X_{i})E(T_{t}-T_{t-1}|\mathscr{F}_{t-1})-\int_{0}^{s}f^{2}(B(v))dv|\xrightarrow[]{P}0.\eqno(5.9)$$
By (5.7) and (5.9), we obtain (5.3).

 {\bf Step 2}
   We prove
     $$\sup_{0<s\le N}|B_{n}(s)-B(s)\circ X_{n}|\xrightarrow[]{P}0, \eqno(5.10)$$
which will be proved if we show
$$\sup_{0<s\le N}|\frac{1}{\sqrt{n}}\sum_{t=2}^{[ns]}f(\frac{1}{\sqrt{n}}\sum_{i=1}^{t-1}X_{i})(X_{t}-D_{t})-\frac{\lambda}{n}\sum_{t=2}^{[ns]}f'(\frac{1}{\sqrt{n}}\sum_{i=1}^{t-1}X_{i})|\xrightarrow[]{P}0.\eqno(5.11)$$

We have
 \begin{eqnarray*}\mathcal
{J}_{k}&:=&|\frac{1}{\sqrt{n}}\sum_{t=2}^{k}f(\frac{1}{\sqrt{n}}\sum_{i=1}^{t-1}X_{i})(X_{t}-D_{t})-\frac{\lambda}{n}\sum_{t=2}^{k}f'(\frac{1}{\sqrt{n}}\sum_{i=1}^{t-1}X_{i})|\\
  &=&|\frac{1}{\sqrt{n}}\sum_{t=2}^{k}f(\frac{1}{\sqrt{n}}\sum_{i=1}^{t-1}X_{i})(H_{t-1}-H_{t})-\frac{\lambda}{n}\sum_{t=2}^{k}f'(\frac{1}{\sqrt{n}}\sum_{i=1}^{t-1}X_{i})|\\
                                     &=& |\frac{1}{\sqrt{n}}\sum_{t=2}^{k}(f(\frac{1}{\sqrt{n}}\sum_{i=1}^{t}X_{i})-f(\frac{1}{\sqrt{n}}\sum_{i=1}^{t-1}X_{i}))H_{t}-\frac{1}{\sqrt{n}}f(\frac{1}{\sqrt{n}}\sum_{i=1}^{k}X_{i})H_{k}-\frac{\lambda}{n}\sum_{t=2}^{k}f'(\frac{1}{\sqrt{n}}\sum_{i=1}^{t-1}X_{i})|,\end{eqnarray*}
and
     \begin{eqnarray*}\max_{1\le k\le nN}\mathcal {J}_{k}&\le & \max_{1\le k\le nN}|\frac{1}{\sqrt{n}}f(\frac{1}{\sqrt{n}}\sum_{i=1}^{k}X_{i})H_{k}|+\max_{1\le k\le nN}|\frac{1}{n}\sum_{t=2}^{k}f'(\frac{1}{\sqrt{n}}\sum_{i=1}^{t-1}X_{i})(X_{t}H_{t}-\lambda)|\\
                                     & & +\max_{1\le k\le nN}|\frac{1}{\sqrt{n}}\sum_{t=2}^{k}(f(\frac{1}{\sqrt{n}}\sum_{i=1}^{t}X_{i})-f(\frac{1}{\sqrt{n}}\sum_{i=1}^{t-1}X_{i})-f'(\frac{1}{\sqrt{n}}\sum_{i=1}^{t-1}X_{i})\frac{X_{t}}{\sqrt{n}})H_{t}|.\end{eqnarray*}
We firstly prove
       $$\max_{1\le k\le nN}|\frac{1}{n}\sum_{t=2}^{k}f'(\frac{1}{\sqrt{n}}\sum_{i=1}^{t-1}X_{i})(X_{t}H_{t}-\lambda)|\xrightarrow[]{P}0.\eqno(5.12)$$
Set $Y_{t,j}=E(X_{t}X_{t+j}|\mathscr{F}_{t})-E(X_{t}X_{t+j})$, we
 prove
                              $$\max_{1\le k\le nN}|\frac{1}{n}\sum_{t=2}^{k}f'(\frac{1}{\sqrt{n}}\sum_{i=1}^{t-1}X_{i})(\sum_{j=1}^{\infty}Y_{t,j})|\xrightarrow[]{P}0.\eqno(5.13)$$
We approximate $S_{t}:=\sum_{j=1}^{\infty}Y_{t,j}$ by
$\widetilde{D}_{t}:=\sum_{k=t}^{\infty}\mathcal {P}_{t}(S_{k})$,
then we need to prove
      $$\max_{1\le k\le nN}|\frac{1}{n}\sum_{t=2}^{k}f'(\frac{1}{\sqrt{n}}\sum_{i=1}^{t-1}X_{i})\widetilde{D}_{t}|\xrightarrow[]{P}0\eqno(5.14)$$
and
      $$\max_{1\le k\le nN}|\frac{1}{n}\sum_{t=2}^{k}f'(\frac{1}{\sqrt{n}}\sum_{i=1}^{t-1}X_{i})(S_{t}-\widetilde{D}_{t})|\xrightarrow[]{P}0.\eqno(5.15)$$
For (5.14), we have
                        $$E(f'(\frac{1}{\sqrt{n}}\sum_{i=1}^{t-1}X_{i})\widetilde{D}_{t})^{2}\le \sqrt{E(f'(\frac{1}{\sqrt{n}}\sum_{i=1}^{t-1}X_{i})^{4})E(\widetilde{D}_{t})^{4}},$$
                         $$[E(\widetilde{D}_{t})^{4}]^{1/4}=[E(\sum_{k=t}^{\infty}\mathcal {P}_{t}(S_{k}))^{4}]^{1/4}\le \sum_{k=t}^{\infty}||{P}_{t}(S_{k})||_{4}.$$
However, by Assumption 3, we have
$$\sum_{k=t}^{\infty}||{P}_{t}(S_{k})||_{4}=\sum_{k=t}^{\infty}\sum_{i=1}^{\infty}||E(X_{k}X_{k+i}|\mathscr{F}_{t})-E(X_{k}X_{k+i}|\mathscr{F}_{t-1})||_{4}<\infty.$$

Since $f'(x)\le C(1+|x|^{\alpha})$,  we have
           $$E(f'(\frac{1}{\sqrt{n}}\sum_{i=1}^{t-1}X_{i})\widetilde{D}_{t})^{2}\le L\eqno(5.16)$$ by Lemma 1. Then, by Kolmogorov inequality for martingale, for any
$\varepsilon>0$ we have
\begin{eqnarray*} & &P(\max_{1\le k\le nN}|\frac{1}{n}\sum_{t=2}^{k}f'(\frac{1}{\sqrt{n}}\sum_{i=1}^{t-1}X_{i})\widetilde{D}_{t}|>\varepsilon)\\
                  &\le &\frac{E[\sum_{t=2}^{nN}f'(\frac{1}{\sqrt{n}}\sum_{i=1}^{t-1}X_{i})\widetilde{D}_{t}]^{2}}{n^{2}\varepsilon^{2}}\\
                    &\le & \frac{N\max_{1\le t\le Nn}E[f'(\frac{1}{\sqrt{n}}\sum_{i=1}^{t-1}X_{i})\widetilde{D}_{t}]^{2}}{n\varepsilon^{2}}\le C\frac{1}{n}\rightarrow 0. \end{eqnarray*}
(5.14) is proved.

For (5.15), we have
             $$S_{t}-\widetilde{D}_{t}=Z_{t-1}-Z_{t},~~Z_{t}=\sum_{i=1}^{\infty}\sum_{k=1}^{\infty}E(X_{t+i+k}X_{t+i}|\mathscr{F}_{t}), \eqno(5.17)$$
and
\begin{eqnarray*}& &\max_{1\le k\le nN}|\frac{1}{n}\sum_{t=2}^{k}f'(\frac{1}{\sqrt{n}}\sum_{i=1}^{t-1}X_{i})(S_{t}-\widetilde{D}_{t})|
                     = \max_{1\le k\le nN}|\frac{1}{n}\sum_{t=2}^{k}f'(\frac{1}{\sqrt{n}}\sum_{i=1}^{t-1}X_{i})(Z_{t}-Z_{t-1})|\\
                     &\le & \max_{1\le k\le nN}|\frac{1}{n}f'(\frac{1}{\sqrt{n}}\sum_{i=1}^{k}X_{i})Z_{k}|+\max_{1\le k\le nN}|\frac{1}{n}\sum_{t=2}^{k-1}(f'(\frac{1}{\sqrt{n}}\sum_{i=1}^{t}X_{i})-f'(\frac{1}{\sqrt{n}}\sum_{i=1}^{t-1}X_{i}))Z_{t}|\\
                      &\le & \max_{1\le k\le nN} |f'(\frac{1}{\sqrt{n}}\sum_{i=1}^{k}X_{i})|\max_{1\le k\le nN} \frac{|Z_{k}|}{n}+\max_{1\le k\le nN} \frac{N|X_{k}Z_{k}|}{\sqrt{n}}\sup_{|t|\le \max_{0\le k\le nN}\frac{\sum_{i=1}^{k}X_{i}}{\sqrt{n}}}f''(t).\end{eqnarray*}
Under Assumption 3, by law of large number, we have
             $$\max_{1\le k\le nN}n^{-\frac{1}{6}}|X_{k}|\xrightarrow[]{P} 0,~~\max_{1\le k\le nN}n^{-\frac{1}{3}}|Z_{k}|\xrightarrow[]{P}0,$$
so $$\max_{1\le k\le nN}
\frac{1}{\sqrt{n}}|X_{k}Z_{k}|\xrightarrow[]{P} 0.$$ From Lemma 2,
we have $\frac{1}{\sqrt{n}}\sum_{i=1}^{k}X_{i}=O_{P}(1)$, by the
continuity of $f''(x)$, we obtain (5.15).

We have, by Talor expansion, that
   \begin{eqnarray*}& &\max_{1\le k\le nN}|\frac{1}{\sqrt{n}}\sum_{t=2}^{k}(f(\frac{1}{\sqrt{n}}\sum_{i=1}^{t}X_{i})-f(\frac{1}{\sqrt{n}}\sum_{i=1}^{t-1}X_{i})-f'(\frac{1}{\sqrt{n}}\sum_{i=1}^{t-1}X_{i})\frac{X_{t}}{\sqrt{n}})H_{t}|\\
                     &\le& (N/2)\max_{1\le k\le nN} \frac{1}{\sqrt{n}}X_{k}^{2}|H_{k}|\sup_{|t|\le \le \max_{0\le k\le nN}\frac{\sum_{i=1}^{k}X_{i}}{\sqrt{n}}}f''(t)\end{eqnarray*}
Under Assumption 3, by law of large number and similar argument in
the above, we get that
$$\max_{1\le k\le nN}|\frac{1}{\sqrt{n}}\sum_{t=2}^{k}(f(\frac{1}{\sqrt{n}}\sum_{i=1}^{t}X_{i})-f(\frac{1}{\sqrt{n}}\sum_{i=1}^{t-1}X_{i})-f'(\frac{1}{\sqrt{n}}\sum_{i=1}^{t-1}X_{i})\frac{X_{t}}{\sqrt{n}})H_{t}|\xrightarrow[]{P}0.$$
Then we can easily get (5.10).

{\bf Step 3}

We  prove (5.2) and $$\lim
_{b\uparrow\infty}\lim\sup_{n}P(|x^{2}|1_{\{|x|>b\}}*\nu_{t\wedge
         S_{n}^{a}}^{n}>\varepsilon)=
         0 \eqno(5.18)$$
          for all $t\in D$, $a>0$.

For (5.2), similar to Ibragimov and Phillips (2008),
    $$\sup_{0<s\le N}|\Delta X_{n}(s)|\le \max_{0\le k\le nN}|f(\frac{1}{\sqrt{n}}\sum_{i=1}^{k}X_{i})|\cdot \max_{0\le k\le nN}\frac{1}{\sqrt{n}}|X_{k}|+\max_{0\le k\le nN}\frac{1}{\sqrt{n}}|X_{k}|.\eqno(5.19)$$
From Lemma 2, we have
$\frac{1}{\sqrt{n}}\sum_{i=1}^{k}X_{i}=O_{P}(1)$. Combining with the
assumptions of $f(x)$, we have (5.2) by (5.19).

As for (5.18), \begin{eqnarray*}& &E\int_{0}^{s\wedge S_{n}^{a}}\int_{\mathbb{R}^{2}}|x^{2}|1_{(|x|>b)}\nu_{n}(dt,dx)\\
                              &\le & E\int_{0}^{s}\int_{\mathbb{R}^{2}}|x^{2}|1_{(|x|>b)}\nu_{n}(dt,dx)\le \frac{1}{b^{2}}E\int_{0}^{s}\int_{\mathbb{R}^{2}}|x^{4}|\nu_{n}(dt,dx)\\
                               &\le &
                               \frac{2}{b^{2}}E\int_{0}^{s}\int_{\mathbb{R}^{2}}|x_{1}^{4}+x_{2}^{4}|\nu_{n}(dt,dx).\\
                           &=&\frac{2}{b^{2}n^{2}}\sum_{t=2}^{[ns]}E[f^{4}(\frac{1}{\sqrt{n}}\sum_{i=1}^{t-1}X_{i})X^{4}_{t}]+\frac{1}{n^{2}}\sum_{t=2}^{[ns]}E[X^{4}_{t}].\end{eqnarray*}

By  uniform continuity of $f(x)$ and
$\frac{1}{\sqrt{n}}\sum_{i=1}^{k}X_{i}=O_{P}(1)$, we have
$f(\frac{1}{\sqrt{n}}\sum_{i=1}^{k}X_{i})=O_{P}(1)$. Furthermore,
$$\frac{1}{n^{2}}\sum_{t=2}^{[ns]}E[f^{4}(\frac{1}{\sqrt{n}}\sum_{i=1}^{t-1}X_{i})X^{4}_{t}]\rightarrow 0$$
and
$$\frac{1}{n^{2}}\sum_{t=2}^{[ns]}E[X^{4}_{t}]\rightarrow 0,$$
we obtain (5.18).

Combining the three steps, we obtain the condition (vi) in Theorem
A, then we get our result.
\section{\bf{Application to Unit Root Autoregression }}
In this section, we use our main result to obtain the limit theorem
for unit root autoregression.  The theory of unit root
autoregression is a hot topic in modern time series.
  Let
   $$Y_{t}=\alpha Y_{t-1}+X_{t},\eqno(6.1)$$
where $X_{t}$ is a causal process with the form
          $$X_{n}=g(\cdots,\varepsilon_{n-1},\varepsilon_{n}).$$
where  $\varepsilon_{n}, n\in Z$ are mean zero, independent and
identically distributed  random variables and $g$ is a measurable
function.

If $\alpha=1$, we want to estimate $\alpha$ from $\{Y_{t}\}$ .
Let$$\hat{\alpha}=\frac{\sum_{t=1}^{n}Y_{t-1}Y_{t}}{\sum_{t=1}^{n}Y^{2}_{t-1}}$$
denote the ordinary least squares estimator of $\alpha$.

Let $t_{\alpha}$ be the regression $t-$statistic with $\alpha=1$:
  $$t_{\alpha}=\frac{(\sum_{t=1}^{n}Y_{t-1}^{2})^{\frac{1}{2}}(\hat{\alpha}-1)}{\sqrt{\frac{1}{n}\sum_{t=1}^{n}(Y_{t}-\hat{\alpha}Y_{t-1})}}.$$
In the following theorem, we get the asymptotic distribution of
$n(\hat{\alpha}-1)$ and $t_{\alpha}$.

{\bf Theorem 2} Under Assumptions 1-3,  we have
  $$n(\hat{\alpha}-1)\xrightarrow[]{d}\frac{\lambda +\sigma^{2}\int_{0}^{1}B(v)dB(v)}{\sigma^{2}\int_{0}^{1}B^{2}(v)dv},\eqno(6.2)$$
  $$t_{\alpha}\xrightarrow[]{d}\frac{\lambda +\sigma^{2}\int_{0}^{1}B(v)dB(v)}{(\int_{0}^{1}B^{2}(v)dv)^{\frac{1}{2}}}.\eqno(6.3)$$

{\bf Proof.}
$$n(\hat{\alpha}-1)=\frac{n\sum_{t=1}^{n}Y_{t-1}X_{t}}{\sum_{t=1}^{n}Y^{2}_{t-1}}.\eqno(6.4)$$
Similarly to  the proof of Theorem 1, we have
    $$(\frac{1}{n}\sum_{t=1}^{[nr]}Y_{t-1}X_{t},\frac{1}{n^{2}}\sum_{t=1}^{[nr]}Y^{2}_{t-1})\Rightarrow(\int_{0}^{r}B^{2}(v)dv,\lambda +\sigma^{2}\int_{0}^{r}B(v)dB(v)).\eqno(6.5)$$
By continuous mapping theorem, we get (6.2) and (6.3).

\section{\bf{Discussion}}

$~~~~~$In this paper, we  study the weak convergence of various
general functionals of partial sums of causal processes.  But we
only consider the univariate case. In Ibragimov and Phillips (2008),
they also considered the multivariate case. However, we can not get
the multivariate extensions by our method. For the univariate case,
we can use the Skorokhod embedding argument to get the asymptotic
results of second predictable characteristics of semimartingale. But
for multivariate case, we can not find a unique stopping time to
embed into  every component of multivariate Brownian motion, so we
can not obtain the correspondence results in accordance with the
method of proof of this paper.

There are mainly two methods to deal with the asymptotic results of
causal processes. One is the martingale approximation developed by
Wu (2007). This method actually is an extension of Beveridge-Nelson
decomposition in Phillips and Solo (1992). The other one is
m-dependent approximation developed by Liu and Lin (2009). By the
m-dependent approximation,  Liu and Lin (2009) get the strong
invariance principle for d-dimensional causal process. The
multivariate extension may be obtained by the m-dependent
approximation in the future research. We take this extension as a
conjecture.

{\bf Conjecture}
    Let $f:\mathbb{R}\rightarrow \mathbb{R}$ be a twice
 continuously differentiable function such that $f'$ satisfying $|f'(x)|\le
 K(1+|x|^{\alpha})$ for some positive constants $K$
and $\alpha$ and all $x\in \mathbb{R}$. Suppose that
$X_{t}=(X^{1}_{t},X^{2}_{t})$ is a 2-dimension  causal process
satisfying  Assumptions 1$\thicksim$3, Then
       $$\frac{1}{\sqrt{n}}\sum_{t=2}^{[nr]}f(\frac{1}{\sqrt{n}}\sum_{i=1}^{t-1}X_{i}^{1})X_{t}^{2}\Rightarrow \lambda\int_{0}^{r}f'(B(v))dv+\sigma\int_{0}^{r}f(B(v))dW(v),$$
where $\lambda=\sum_{j=1}^{\infty}E|X_{0}X_{j}|$, $(B(s), W(s))$ is
a bivariate Brownian motion.
\section*{Reference}
{\small Billingsley P. (1968). {\em Convergence of Probability
Measure.}
Wiley.\\
Caceres C, Nielsen B. (2007) Convergence to Stochastic Integrals
with Non-linear Integrands. Manuscrips.
\\Chan N-H, Wei C-Z. (1987). Asymptotic
inference for nearly nonstationary AR(1) processes. {\em
Annals of Statistics 15}, 1050-1063.\\
Chan N-H, Wei C-Z. (1988). Limiting distribution of least square
estimates of unstable autoregression processes.{\em
Annals of Statistics 16}, 367-401.\\
De Jong R, Davidson J. (2000 a). The functional central limit
theorem and weak convergence to stochastic integral I: weak
dependent
processes. {\em Econometric Theory 16}, 621-642.\\
De Jong R, Davidson J. (2000 b). The functional central limit
theorem
and weak convergence to stochastic integral II: fractionally integrated processes.  {\em Econometric Theory 16}, 643-666.\\
Dellacherie  C, P Meyer.  (1982)  {\em Probability and Potential B.}
North-Holland.\\
 He S-W, Wang J-G, Yan J-A. (1992). {\em  Semimartingale and
Stochastic
Calculus.} CRC Press.\\
Ibragimov R, Phillips P. (2008). Regression asymptotics using
martingale convergence methods. {\em Econometric Theory 24},
888-947.\\
Ikeda J. Watanabe S. (1989). {\em Stochastic Differential Equations
and Diffusion Processes. } North-Holland.\\
 Jacod J, Shiryaev AN.
(2003). {\em Limit Theorems for Stochastic
   Processes. } Springer. \\
Jeganathan P. (2004). Convergence of functionals of sums of r.v.s to
local times of fractional motions. {\em
Annals of Probability 32}, 1771-1795.\\
Liu W-D, Lin Z-Y. (2009). Strong approximation for a class of
stationary processes. {\em Stochastic Processes and their
Applications
 119}, 249-280. \\
Phillips P.C.B. (1987 a). Time-series regression with a unit root.
{\em Econometrica 55}, 277-301.\\
Phillips P.C.B. (1987 b).  Towards a unified asymptotic theory for
autoregression. {\em Biometrika 74}, 535-547.\\
Phillips P.C.B. (2007). Unit root log periodogram regression. {\em
Journal of Econometrics 138}, 104-124.\\
Phillips P.C.B, Solo V. (1992). Asymptotic for linear process. {\em
Annals of Statistics 20}, 971-1001.\\
Wu, W-B. (2005) Nonlinear system theory: Another look at
 dependence. {\em Proceedings of the National Academy of Science
 102}, 14150-14154.\\
 Wu W-B.(2007). Strong invariance
principles for dependent random variables. {\em The Annals of
Probability  35}, 2294-2320. }

 \end{document}